\documentclass[1p,number]{elsarticle}
\usepackage{amsmath}
\usepackage{amssymb}
\usepackage{amsfonts}
\usepackage{amsthm}
\usepackage{latexsym, ulem}
\usepackage{graphicx}
\begin{document}

\begin{frontmatter}
\title{A new identity linking coefficients of a certain class of homogeneous polynomials and Gauss hypergeometric functions}
\author{Philip W. Livermore}
\ead{phil@ucsd.edu}
\author{Glenn R. Ierley}
\address{Institute of Geophysics and Planetary Physics, Scripps Institution of Oceanography, UCSD, La Jolla, CA 92093 USA}

%    Abstract is required.
\begin{abstract}
In investigating the properties of a certain class of homogeneous polynomials, we discovered an identity satisfied by their coefficients which involves simple ${}_2F_1$ Gauss hypergeometric functions. This result appears to be new and we supply a direct proof.
The simplicity of the identity is suggestive of a deeper result.
\end{abstract}

\end{frontmatter}

The results that we report here arose when investigating an applied mathematics problem ostensibly unrelated to combinatorics:  the self-generation of magnetic fields in a spherical geometry as applicable, for instance, to the Earth's fluid core.
It happens that a certain set of homogeneous conditions apply to such magnetic fields, the so-called Taylor constraints, which arise in the physically interesting low-viscosity fast-rotation limit of the governing equations \citep{Taylor_63}. These constraints take the form of the vanishing of certain integrals over surfaces of constant cylindrical radius of a quantity involving the magnetic field. We have recently proven that a necessary condition for any solution is the vanishing of a certain function of cylindrical radius on both inner and outer spherical boundaries of the core \cite{Livermore_etal_2009_PEPI_submitted}.
On adopting a truncated spatial discretisation based on spherical harmonics in solid angle and certain regular polynomials in radius, this amounts to the vanishing of
\begin{equation}  {\cal Q}(\rho,s) = \sum_{j,k} a_{jk} \, {\cal I}_{jk}(\rho,s)  \label{eqn:overall} \end{equation}
regarded as a function of $s$ at fixed values of $\rho=7/20$ and $\rho=1$. The coefficients $a_{jk}$ are real and 
\begin{equation} {\mathcal I}_{jk}(\rho,s) = \int_0^{\sqrt{\rho^2-s^2}} z^{2j} (z^2+s^2)^k\, dz. \label{eqn:form_fn} \end{equation}
In \eqref{eqn:overall}, $j$ and $k$ belong to a prescribed finite set (whose size depends on the spectral truncation adopted in the numerical scheme) and $s$ is the nondimensional cylindrical radius; the given values of $\rho$ represent respectively the spherical radii of the inner and outer boundaries of the Earth's fluid core. 
We note two points associated with any given $(j,k)$. Firstly, the integrand in \eqref{eqn:form_fn} is a homogeneous polynomial in $(s,z)$ of degree $2(j+k)$; secondly, its limits are homogeneous of degree one. Thus ${\mathcal I}_{jk}(\rho,s)$ is homogeneous of degree $2(j+k)+1$ and can be written
\begin{equation}
{\mathcal I}_{jk}(\rho,s) = \sqrt{\rho^2-s^2}\sum_{l=0}^{j+k} B_l(\rho)\, s^{2l}. \label{eqn:Proof1} \end{equation}
It follows immediately that, for each choice of $\rho$, \eqref{eqn:overall} becomes 
\begin{equation} {\cal Q}(\rho,s) = \sqrt{\rho^2-s^2} \sum_{l=0}^L \, A_l(\rho) \, s^{2l}  \label{eqn:overall2}\end{equation}
where $L$ is the maximum value of $j+k$.

We may therefore impose the required conditions, the vanishing of both ${\cal Q}(7/20,s)$ and ${\cal Q}(1,s)$, simply by setting each coefficient $A_l$ appearing in \eqref{eqn:overall2} individually to zero. This procedure yields a set of constraints of size $2L$. 
However, it happens that this set contains (in general) many degeneracies. 
Empirically, we found that the $A_l$, produced by choosing one of the two values of $\rho$ as unity, are related by
\begin{equation} A_N(1) = A_N(\rho) + \frac{1-\rho^2}{2} \sum_{m=1}^{L-N} \mu_m(\rho)\, A_{N+m}(\rho),\qquad N \ge \max \{j\} \label{eqn:empirical1} \end{equation}
where $\mu_m$ is proportional to a ${}_2F_1$ Gauss hypergeometric, given below.   In our Taylor-state problem, we will fix $\rho=7/20$, but the above statement is true more generally.
Thus, for any $N \ge \max\{j\}$, if the values of $A_l$ appearing on the right hand side above are set to zero, then $A_N(1)$ automatically vanishes and explicitly enforcing $A_N(1) = 0$ is unnecessary and introduces a degeneracy in the set of constraints.

The source of the linear homogeneous condition \eqref{eqn:empirical1} is immediately traced to the same condition on each ${\cal I}_{jk}$ individually, that is,
\begin{equation} B_N(1) = B_N(\rho) + \frac{1-\rho^2}{2} \sum_{m=1}^{j+k-N} \mu_m(\rho)\, B_{N+m}(\rho), \qquad N \ge j \label{eqn:empirical2} \end{equation}
where $B_l$, defined in \eqref{eqn:Proof1}, is proportional to $\rho^{2(j+k-l)}$ as follows from homogeneity. Note that the upper limit on the summation is the tightest possible since $B_l = 0$ if $l > k+j$.
It is immediate that \eqref{eqn:empirical1} follows from \eqref{eqn:empirical2} because the summand is independent of $(j,k)$. 

In this note we provide a proof of a generalization of \eqref{eqn:empirical2}, and it is then a simple matter to use \eqref{eqn:empirical1} to count the number of independent conditions related to enforcing Taylor's constraint.

The particular structure of the integrand arises in the following manner.
The monomial dependence of $z$ stems from the alignment of the cylindrical integrals, defining Taylor's constraint, with the $z$-axis. The appearance of $s^2$ only through $s^2+z^2$ comes about due to assumed $C^{\infty}$ behaviour of the magnetic field: since both $z$ and the square of the spherical radial distance, $r^2 = s^2 +z^2$, are both $C^{\infty}$, any functional relation of the given form also inherits this property. Lastly, the form of the upper limit on integration comes from the height of the intersection of a cylinder of cylindrical radius $s$ with a sphere of radius $\rho$.
The elementary structure of the identity is suggestive of a deeper result, although the statement cannot be readily generalized, for instance, to odd exponents of $z$.  It is possible that this result arises in other circumstances involving cylindrically symmetric constraints in fluid dynamics, for instance, the Taylor-Proudman condition in a rotating sphere \cite{Book_Batchelor_67}. 

Immediately below we state the main result, but defer proof until two necessary Lemmas are stated and proved.
\newtheorem{thm}{Theorem}
\begin{thm}
Let $j$ and $k$ be positive integers and let the coefficients $B_l(\rho)$ be defined by
\begin{equation} {\mathcal I}_{jk}(\rho,s) = \int_0^{\sqrt{\rho^2-s^2}} z^{2j} (z^2+s^2)^k\, dz = \sqrt{\rho^2-s^2}\sum_{l=0}^{j+k} B_l(\rho)\, s^{2l}. \end{equation}

Then, for any positive integer $N$ with $k+j \ge N \ge j$,  the quantity 
\begin{equation} B_N(\rho) + \frac{1-\rho^2}{2} \sum_{m=1}^{j+k-N} \mu_m(\rho)\, B_{N+m}(\rho) \label{eqn:rho_indep2} \end{equation}
is independent of $\rho$ where
\[ \mu_m(\rho) =  \frac{2\Gamma(m+1/2)} {\sqrt{\pi}\, \Gamma(m+1)} {}_2F_1\left ( \genfrac{}{}{0pt}{}{[1-m,1/2]}{[1/2-m]};\rho^2 \right).\]
\end{thm}

\newdefinition{rmk1}{Remarks}
\begin{rmk1}
The statement in the theorem is considerably stronger than that of equation \eqref{eqn:empirical2} and relates expressions of the form \eqref{eqn:form_fn} between any two values of $\rho$, not necessarily including unity.
Additionally, since $\mu_0 = 2 (1-\rho^2)^{-1}$, the theorem is equivalent to the $\rho$-independence of
 \[ \frac{1-\rho^2}{2}\; \sum_{m=0}^{j+k-N} \mu_m B_{N+m}.\]
Although this is the most succinct form of the result, the statement in the theorem is easier to prove as $\mu_m(\rho)$, for $m \ge 1$, are simply polynomials.
The $\mu_m$ have been arbitrarily normalised such that $\mu_0 = 1$.
\end{rmk1}

\newtheorem{lem}[thm]{Lemma}
\begin{lem}
For any real values $\alpha,\beta$, the following are identities
\begin{align} 
\sum_{k=0}^n (-1)^{n+k} \binom{\alpha+k}{\beta+n}\, \binom{n}{k} &= \binom{\alpha}{\beta}, \label{eqn:lem1}  \\
\sum_{k=0}^n \frac{1}{\beta+k} \binom{n-k-1/2}{-1/2}\, \binom{k-1/2}{-1/2} & ={\frac {\Gamma  \left( \beta \right) \Gamma  \left( n+\beta+1/2  \right) }{\Gamma  \left( \beta+1/2 \right) \Gamma  \left( \beta+n+1  \right) }},  \label{eqn:lem2} 
%\sum_{k=0}^n \frac{ \Gamma(n-k+1/2)\, \Gamma(\alpha+k)\,\Gamma(k+1/2)}{\Gamma(n-k+1)\,\Gamma(\alpha+k+1)\,\Gamma(k+1)} = \frac{\pi\,\Gamma(\alpha)\,\Gamma(n+\alpha+1/2)}{\Gamma(\alpha+1/2)\,\Gamma(n+\alpha+1)}
\end{align}
\end{lem}
where $\binom{\alpha}{\beta}$ is a binomial coefficient (extended to non-integer values of its arguments).

\newproof{pflem}{Proof}
\begin{proof}
We rewrite the left hand side of the above two equations in a hypergeometric form,
\begin{align}
&\frac{ (-1)^n\Gamma(\alpha+1)}{\Gamma(\alpha-\beta-n+1)\,\Gamma(\beta+n+1)}\, {}_2F_1\left ( \genfrac{}{}{0pt}{}{[-n,\alpha+1]}{[\alpha-\beta-n+1]};1 \right), \\
&\frac{\Gamma(n+1/2)}{\sqrt{\pi}\,\Gamma(n+1)\,\beta} \, {}_3F_2\left( \genfrac{}{}{0pt}{}{[-n,1/2,\beta]}{[-n+1/2,\beta+1]}  ;1\right).
\end{align}
Using the theorems of Gauss and  Saalsch\"{u}tz \cite{Book_Petkovsek_etal_97}, we may evaluate the ${}_2F_1$ and ${}_3F_2$ functions and the results follow easily.
\end{proof}

\newproof{pf}{Proof of main result}
\begin{pf} 

We prove the theorem directly in two steps. First, we find a general expression for the coefficients $B_l$ and second, we show the required relation holds between them.
\begin{enumerate}
\item[Step 1]
By exploiting equation B6 in \cite{Livermore_etal_2008}, we can write
\[ {\mathcal I}_{jk}(\rho,s) = \frac{ \rho^{2k} }{2j+2k+1}\, \left(\rho^2-s^2\right)^{j+1/2} \, {}_2F_1\left ( \genfrac{}{}{0pt}{}{[1,-k]}{[1/2-j-k]}; s^2/\rho^2 \right).\]
The hypergeometric function appearing is simply a polynomial in $s^2/\rho^2$ of degree $k$, 
\[  {}_2F_1\left ( \genfrac{}{}{0pt}{}{[1,-k]}{[1/2-j-k]}; s^2/\rho^2 \right) = \sum_{n=0}^k \frac{\Gamma(k+1) \, \Gamma(k+j-n+1/2)}{\Gamma(k-n+1) \, \Gamma(k+j+1/2)} \left(\frac{s}{\rho}\right)^{2n}\] 

and ${\mathcal I}_{jk}(\rho,s)$ can be written
\begin{equation} {\mathcal I}_{jk}(\rho,s)  =  \left( \rho^2-s^2 \right)^{1/2} \frac{\Gamma(k+1)}{2\,\Gamma(k+j+3/2)}  \left( \rho^2-s^2 \right)^{j}\,\sum_{n=0}^k \frac{ \Gamma(k+j-n+1/2)}{\Gamma(k-n+1)} s^{2n}\, \rho^{2(k-n)}  \nonumber \end{equation}
where, aside from the leading prefactors,  the binomial and summation appearing immediately above can be further written as
\begin{equation} \sum_{n=0}^k \,\sum_{i=0}^j \frac{ \Gamma(k+j-n+1/2)\,\Gamma(j+1)\,  (-1)^i}{\Gamma(k-n+1)\,\Gamma(j-i+1)\, \Gamma(i+1)} s^{2(n+i)}\, \rho^{2(k-n+j-i)}. \end{equation}
Up to the prefactor of $ \left( \rho^2-s^2 \right)^{1/2}$, ${\mathcal I}_{jk}(\rho,s)$ is simply a polynomial in $s^2$ of degree $k+j$. To express this in the form of \eqref{eqn:Proof1}, we may write $l=n+i$ and re-order the sum over dummy indices $l$ and $n$.
However, care must be taken with the limits. Noting that $i \ge 0$ and therefore $n \le l$, coupled with the upper bound $n\le k$, leads to $n \le \min(k,l)$. Additionally, since $i \le j$ and therefore $n\ge l-j$, coupled with the lower bound $n \ge 0$, leads to $n \ge \max(0,l-j)$. Hence
\begin{equation}  B_l(\rho) = \frac{ \Gamma(k+1)\,\Gamma(j+1)}{2\,\Gamma(k+j+3/2)} \sum_{n=\max(0,l-j)}^{\min(k,l)} \frac{ (-1)^{l+n} \, \Gamma(k+j-n+1/2) }{\Gamma(k-n+1)\,\Gamma(j+n-l+1)\,\Gamma(l-n+1)}\rho^{2(k+j-l)}. \label{eqn:Bl_eqn} \end{equation}
Shifting the dummy indices by writing $m=n-(l-j)$ and defining $T=k+j-l$, this expression can be simplified to
\begin{equation}
B_l(\rho) = \frac{ \Gamma(k+1)\,\Gamma(j+1/2)\,\rho^{2(k+j-l)}}{2\,\Gamma(k+j+3/2)} \sum_{m=\max(0,j-l)}^{\min(T,j)} (-1)^{m+j}\,\binom{T-1/2+j-m}{j-1/2} \binom{j}{j-m}. \label{eqn:Bl_explicit} \end{equation}
Note that if $m>T$ or $m>j$ then respectively, the first or second binomial term vanishes. Hence, without loss of generality, we can replace the upper limit of the sum by $j$. If $l \ge j$ then the lower limit is zero and a trivial reordering of the sum, setting $m' = j-m$, and using \eqref{eqn:lem1}, leads to
\begin{equation}
B_l(\rho) = \frac{ \Gamma(k+j-l+1/2)\, \Gamma(k+1)\, \Gamma(j+1/2)\,  (-1)^j}{2\sqrt{\pi} \,\Gamma(k+j+1-l) \, \Gamma(k+j+3/2)} \rho^{2(k+j-l)}. 
\end{equation}
If $l < j$ then the summation in \eqref{eqn:Bl_explicit} has a nonzero lower limit and cannot be evaluated in such a simple closed form.

% main result.
\item[Step 2]
We now are in a position to prove that, if $l \ge j$ then the $B_l(\rho)$ satisfy \eqref{eqn:rho_indep2}.
We shall show that
\begin{equation} \sum_{m=1}^{k+j-N} B_{N+m}(\rho) \,\mu_m(\rho) = \frac{ \Gamma(k+j-N+1/2)\, \Gamma(k+1)\, \Gamma(j+1/2)}{\sqrt{\pi} \,\Gamma(k+j+1-N) \, \Gamma(k+j+3/2)} \sum_{b=0}^{k+j-N-1} \rho^{2b} \label{eqn:required1} \end{equation}
and it follows immediately that
\begin{align}
&B_N(\rho) + \frac{(1-\rho^2)}{2} \, \sum_{m=1}^{k+j-N} B_{N+m}(\rho) \,\mu_m(\rho) \nonumber\\
&=  \frac{ \Gamma(k+j-N+1/2)\, \Gamma(k+1)\, \Gamma(j+1/2)}{2\sqrt{\pi} \,\Gamma(k+j+1-N) \, \Gamma(k+j+3/2)} \left(  \rho^{2(k+j-N)} + (1-\rho^2) \sum_{b=0}^{k+j-N-1} \rho^{2b} \right) \nonumber \\
& = \frac{ \Gamma(k+j-N+1/2)\, \Gamma(k+1)\, \Gamma(j+1/2)}{2\sqrt{\pi} \,\Gamma(k+j+1-N) \, \Gamma(k+j+3/2)}, \end{align}
independent of $\rho$.
Note that \eqref{eqn:rho_indep2} involves $B_l$ with $l \ge N$. Coupled with the initial hypothesis $N \ge j$, this is consistent with $l \ge j$, required in the derivation of the closed form for $B_l$.

%The origin of the bound $N \ge j$ stems from the fact that only $B_N, B_{N+1},\dots$ arise in \eqref{eqn:rho_indep2} and the required simplication of the summation in \eqref{eqn:Bl_eqn} is only possible when the indexes are at least $j$, i.e. $N \ge j$. 

It remains to show \eqref{eqn:required1}. By expanding the hypergeometric component of $\mu_m(\rho)$, a polynomial of degree $m-1$ in $\rho$, it follows easily that 
\[ \mu_m(\rho) =  \sum_{n=0}^{m-1}  \frac{ 2 }{\pi} \frac{\Gamma(m)\, \Gamma(n+1/2)\, \Gamma(m-n+1/2)}{ \Gamma(m+1)\,\Gamma(m-n)\, \Gamma(n+1)}\rho^{2n} .\]
Hence
\begin{align} &\sum_{m=1}^{k+j-N} B_{N+m}(\rho) \,\mu_m(\rho) = \frac{(-1)^j\,\Gamma(k+1)\, \Gamma(j+1/2)}{\pi^{3/2}\, \Gamma(k+j+3/2)} \times \nonumber\\
& \sum_{m=1}^{k+j-N} \;  \sum_{n=0}^{m-1} \frac{ \Gamma(k+j-N-m+1/2)\, \Gamma(m)\, \Gamma(n+1/2)\, \Gamma(m-n+1/2) }{\Gamma(k+j+1-N-m) \, \Gamma(m+1)\,\Gamma(m-n)\, \Gamma(n+1)}  \rho^{2(k+j+n-N-m)}. \label{eqn:long_formula} \end{align}
This is a polynomial in $s^2$ of degree $k+j-N-1$.
By introducing a new dummy variable $b=k+j+n-N-m$ the double sum above can be rewritten as
\[ \sum_{m=1}^{k+j-N} B_{N+m}(\rho) \,\mu_m(\rho) = \sum_{b=0}^{k+j-N-1} \nu_b\, \rho^{2b} \]
where 
\begin{align} &\nu_b =\frac{(-1)^j\,\Gamma(k+1)\, \Gamma(j+1/2)\,\Gamma(k+j-N-b+1/2)}{\pi^{3/2}\, \Gamma(k+j+3/2)\,\Gamma(k+j-N-b)} \times \nonumber \\
& \sum_{m=k+j-N-b}^{m=k+j-N} \frac{ \Gamma(k+j-N-m+1/2)\, \Gamma(m)\, \Gamma(b-k-j+N+m+1/2) }{\Gamma(k+j+1-N-m) \, \Gamma(m+1)\, \Gamma(b-k-j+N+m+1)}. \label{eqn:required_result2}
\end{align}
The change in the lower limits on $m$ arises since $b-k-j+m+N = n \ge 0$ in the original sum and hence $m \ge k+j-N-b$. Note that, since $b \le k+j-N-1$, $m \ge 1$, consistent with the original lower limit on $m$ from \eqref{eqn:long_formula}.
Defining new variables $c=m-\beta$, $\beta =  k+j-N-b$, the summation can be written
\[ \sum_{c=0}^b \frac{\Gamma(b+c+1/2)\,\Gamma(c+1/2)}{\Gamma(b+c+1)\,\Gamma(c+1)}\frac{1}{c+\beta} = \sum_{c=0}^b \binom{b+c-1/2}{-1/2}\,\binom{c-1/2}{-1/2} \frac{\Gamma(1/2)^2}{c+\beta}  \]
and by using \eqref{eqn:lem2} it is immediate that
\[ \nu_b =   (-1)^j\,\frac{ \Gamma(k+j-N+1/2)\, \Gamma(k+1)\, \Gamma(j+1/2)}{\sqrt{\pi} \,\Gamma(k+j+1-N) \, \Gamma(k+j+3/2)}\]
which is independent of $b$, and \eqref{eqn:required1} follows.
%\[ \frac{ \Gamma(M-N+1/2)}{\Gamma(M-N+1)}\,\sum_{e=0}^{W} \frac{\Gamma(W-e+1/2)\,\Gamma(e+1/2)}{\Gamma(W-e+1)\,\Gamma(e+1)\,(e+M)}\]
%\[\frac{ \Gamma(M-N+1/2)}{\Gamma(M-N+1)}\,\sum_{e=0}^{W} \binom{W-e-1/2}{-1/2}\,\binom{e-1/2}{-1/2} \frac{\Gamma(1/2)^2}{(e+M)}\]
%
\end{enumerate}
\end{pf}

\baselineskip=11pt

\providecommand{\bysame}{\leavevmode\hbox to3em{\hrulefill}\thinspace}
\providecommand{\MR}{\relax\ifhmode\unskip\space\fi MR }
\providecommand{\MRhref}[2]{%
  \href{http://www.ams.org/mathscinet-getitem?mr=#1}{#2}
}
\providecommand{\href}[2]{#2}

\end{document}